\documentclass[12pt]{article}
\usepackage{amsmath,amstext, amsfonts, verbatim,remreset}
\usepackage{graphicx}

\newtheorem {theorem} {Theorem}
\newtheorem {lemma} {Lemma}

\newcommand {\RR} {\mathbb R}
\newcommand {\NN} {\mathbb N}

\newcommand {\ZZ} {\mathbb Z}

\newcommand {\xx} {\bf x}

\font\mathcal=cmbsy10 

\begin{document}

\title{Exact asymptotics of the uniform error of interpolation by multilinear splines}


\author{Yuliya Babenko\\ Sam Houston State University}
\date{}
\maketitle

{{ \centerline{\bf Abstract} The question of adaptive mesh generation for approximation by splines has been
studied for a number of years by various authors. The results have numerous applications in computational and
discrete geometry, computer aided geometric design, finite element methods for numerical solutions of partial
differential equations, image processing, and mesh generation for computer graphics, among others.
In this paper we will investigate the questions of adaptive approximation of $C^2$ functions with arbitrary but fixed throughout the domain signature by multilinear splines. In particular, we
will study the asymptotic behavior of the optimal error of the weighted uniform approximation by interpolating and quasi interpolating multilinear splines. }}


\section{Introduction.}

Let the domain $D$ be the unit cube $[0,1]^d \subset \RR^d$, $d \in \NN$ and $d\geq 2$. However, any bounded
connected region that can be represented as a finite union of cubes can be treated analogously.

Let $L_{\infty}(D)$ be the standard space of essentially bounded measurable functions defined on $D$ with the
usual sup-norm $\|\cdot\|_{\infty}$. Given a positive continuous function $\Omega({\bf x})$
 on $D$, define the weighted uniform
norm $\|\cdot\|_{\infty, \Omega}$ of $f \in L_{\infty}(D)$ as
$$
\|f\|_{\infty,\Omega}:= \displaystyle {\rm{{\rm \sup}}}\left\{\; |f({\bf x})|\Omega({\bf x})\; : \; {{\bf x} \in
D}\;\right\}=\|f\Omega\|_{\infty}.
$$

 Let $\square_N=\{R_i \}_{i=1}^N$ be an arbitrary partition of the domain $D$ into $N$ parallelepipeds (or boxes) $R_i$ with sides parallel to coordinate axes. We shall refer to such a partition as a {\it box partition. }

Let 
$
{\cal {P}}^*_1
$
 be the space of polynomials linear in each of its $d$ variables.
Given a function $f \in L_{\infty}(D)$ and a partition $\square_N$ we shall consider the multilinear spline $  s (f, \square_N)$ defined as follows:
\begin{enumerate}
\item On the interior of each box $R_i,\; i=1, \dots, N$, we define  $  s (f, \square_N)$ to be the polynomial from ${\cal {P}}^*_1$ which interpolates $f$ at the vertices of $R_i$.

\item For every point  $x$  in the union of the boundaries of boxes $R_i$, $i=1,\dots, N$ we define the value of $s (f, \square_N)$ to be the average value at $x$ of all polynomial interpolants on the boxes whose boundaries contain $x$.

\end{enumerate}

Observe that, by construction, splines $  s(f, \square_N)$ interpolate the function $f$ at the vertices of the partition $\square_N$. However, in general, they are not necessarily continuous. 

We shall call a sequence of partitions $\{\square_N\}_{N=1}^{\infty}$ {\it admissible} if it satisfies the condition
\begin{equation}\label{asss}
\sup_N  N^{\frac 1d} \displaystyle \max _{R \in \square_N} {\rm diam}(R) < \infty.
\end{equation}
Throughout this paper we shall consider only admissible box partitions.

The quantity
\begin{equation} \label{2}
{\mathbf R}_N(f):=\inf_{\square_N} \|f-  s(f,\square_N)\|_{\infty,\Omega}
\end{equation}
where $\inf$ is taken over all admissible box partitions $\square_N$ of the domain $D$ into $N$ boxes
we shall call the {\it optimal error} of interpolation.
An explicit form and the exact value of ${\mathbf R}_N(f)$, as well as the explicit construction of the optimal
partition for every particular function $f$, can be found only in trivial situations. It was shown by Below, De Loera,
and Richter-Gebert in 2000 ~\cite{DeL} that it is not possible to construct an adaptive algorithm for optimal mesh
(partition) generation that runs in polynomial time.

That is why it is interesting to study the asymptotics of the optimal error ${\mathbf R}_N(f)$
 as $N \to \infty$  for a given function $f \in C^2$ and to construct an asymptotically optimal sequence of box partitions, i. e. a
 sequence of partitions $\{ \square_N^* \}^{\infty}_{N=1}$ of $D$ such that
\begin{equation}\label{3}
\displaystyle \lim_{N \to \infty}\frac{\|f-  s(f,\square^*_N)\|_{\infty,\Omega} }{{\mathbf R}_N(f)}=1.
\end{equation}

Note that the problem formulated above is interesting for functions of arbitrary smoothness as well as for various
classes of splines (for instance, for splines of higher order, interpolating splines, best approximating splines, best
one-sided approximating splines, etc.). In the univariate case general questions of this type have been investigated by
many authors. The results obtained in this case are more or less complete and have numerous applications (see, for
example, ~\cite{LSh}).

 Fewer results are known in the multivariate case. The
classical statement of L. Fejes Toth indicated in ~\cite{Toth} on approximation of convex bodies by inscribed polytopes
in Hausdorff metric can be considered as the first result in this direction. Gruber ~\cite{Gr2} proved this result and
generalized it to the arbitrary dimension.  We ~\cite{us} proved similar result for the weighted uniform norm. Related interesting results on approximation of convex bodies by various
polytopes have been obtained by B$\rm{\ddot{o}}$r$\rm{\ddot{o}}$czky ~\cite{Bor}, B$\rm{\ddot{o}}$r$\rm{\ddot{o}}$czky
and Ludwig ~\cite{KL}.

In ~\cite{Nadler} Nadler  solved the problem of asymptotically optimal choice of a sequence of triangulations for
approximation of $C^3$ functions by piecewise linear functions of best $L_2$-approximation.

D'Azevedo and Simpson ~\cite{Daz3} studied the question of triangulating a given set of vertices for interpolation of a
convex quadratic surface by piecewise linear functions. They showed that the Delaunay triangulation will be optimal for
the error in the $L_{\infty}$ norm. For the error in $L_p$ norm this fact was proved by Rippa ~\cite{Rippa}. Chen, Sun, and
Xu ~\cite{Chen1} generalized this result to arbitrary dimensions.

Later D'Azevedo ~\cite{Daz2} obtained local error estimates for functions with both positive and negative curvature.
The same estimates were later obtained by Pottmann, Hamann {\it {et al}} ~\cite{kodla1} who studied the problem of
optimally triangulating the plane for approximating quadratic functions by piecewise linear functions and suggested
some algorithms for constructing function dependent triangulations of the whole domain.

Huang ~\cite{Huang}, and Huang, Xu, and Sun ~\cite{HSun} considered the problem of variational mesh adaptation in the
numerical solutions of partial differential equations and  obtained asymptotic bounds on the interpolation error
estimates in $L_2$ for adaptive meshes that satisfy regularity and equidistribution conditions.

All above mentioned results are for the case of approximation by linear splines. Natural domain partitions in this case are
simplices. However, in applications where preferred directions exist, box partitions (or generalized rectangular
partitions) are sometimes more convenient. The only known to us result for box partitions is the estimate of the uniform error of
interpolation on rectangular partition by bivariate splines linear in each variable which is due to D'Azevedo
~\cite{Daz1} who obtained the local error estimates in this case.

Therefore, we think it is an interesting problem to obtain results similar to the above mentioned for the class of splines linear in each variable (multilinear splines) defined over box partitions.

 In this paper we shall construct an asymptotically optimal sequence $\{\square^*_N\}_{N=1}^{\infty}$ of box partitions and determine the exact asymptotic value of ${\mathbf R}_N(f)$ for a function $f \in C^2(D)$ with fixed signature (with fixed number of positive and negative second derivatives) throughout the domain. Even though signature is fixed, the results presented are richer than what has been proved for interpolating by linear splines on $\RR^d$, $d>2$. In the latter case, while some results exist for positive definite functions, no results exist in the case of functions of another signature. 
 
 The approach we shall undertake is similar to the one used for studying the asymptotics of the error of approximation by linear splines and is as follows: we first take an intermediate approximation of the given function $f$ by a quadratic polynomial and then find the error of approximating the quadratic by multilinear splines. This last problem is solved for quadratic functions of arbitrary fixed signature. The result is then used to give the error of approximating a $C^2$ function with (fixed throughout the domain) arbitrary signature. Even though it has not being done in the current text, it is rather clear how to proceed and extend the obtained results to approximate an arbitrary $C^2$ function.

In addition, we show that the corresponding multilinear splines $ \{ s(f, \square^*_N)\}_{N=1}^{\infty}$ can be constructed so that they will be discontinuous only along small number (in comparison with the total number) of lines. 

If we do not require interpolation at every vertex of a partition we can construct an asymptotically optimal sequence of admissible partitions $\{\tilde \square^*_N\}_{N=1}^{\infty}$ and a sequence of continuous splines  $ \{ \tilde s(f, \tilde \square^*_N)\}_{N=1}^{\infty}$ which interpolate $f$ at all but $o(N)$ vertices of the partition as $N \to \infty$. We shall refer to such splines as {\it quasi interpolating} splines.

\section{ Main results and ideas of proofs.}

For $f \in C^2(D)$ denote
\begin{equation}
H(f;{\bf x}):=\displaystyle \prod_{i=1}^d \frac{\partial ^2 f}{\partial x_i^2}({\bf x}).
\end{equation}
Let us consider for any $0\leq k \leq d$ the following class of functions
$$
C^2_k(D):=\left \{ f\in C^2(D) : \forall {\bf x}\in D \;\frac {\partial ^2 f}{\partial x_i^2}({\bf x})>0, \; 1\leq i \leq k, \;\hbox{and}\;\frac{\partial ^2 f}{\partial x_i^2}({\bf x})<0,\; k < i \leq d\; \right \}.
$$
Sometimes, we shall say that functions from $C^2_k(D)$ have signature $(k,d)$. In the case when $k=0$ or $k=d$ we shall say that functions are positive definite.

In addition, for $k, d\in \NN \cup \{0\}$ set
\begin{equation} \label{gamma}
\gamma_{k,d}:=\begin{cases}
\frac 18 k^{\frac
kd}(d-k)^{1-\frac kd}, & 0<k<d \cr
\frac{d}{8}, & k=d \;\hbox{or}\; k=0. \cr
\end{cases}
\end{equation}

The next theorem contains the main result of this paper.

\begin{theorem} \label{L} 
For any $0\leq k \leq d$ and $f \in C_k^2(D)$
\begin{equation} \label{d_tthm_par_b}
\displaystyle  \lim_{N \to \infty} N^{\frac 2d}{\mathbf R}_N(f) = \gamma_{k,d}\left ( \displaystyle \int_{D} |H(f;{\bf x})|^{\frac{1}{2}}\Omega({\bf {x}})^{\frac d2} d{\bf
x}\right)^{\frac{2}{d}}.
\end{equation}

Furthermore, there exists a sequence of admissible box partitions $\{\tilde  \square_N\}_{N=1}^{\infty}$ and a sequence of continuous quasi interpolating splines $\tilde s(f, \tilde \square_N)$ such that
$$
\displaystyle \lim_{N \to \infty} \frac{\|f-\tilde s(f,\tilde \square_N)\|_{\infty, \Omega}}{{\mathbf R}_N(f)}=1.
$$
\end{theorem}

{\bf Remark.} The sequence of splines constructed in the proof of the upper bound in (\ref{d_tthm_par_b}) will possess the following nice property: for each $N$ the constructed spline will be discontinuous along only a small number,  i.e. $o(N)$ of faces as $N\to \infty$ (see Section 5.2).

Let us describe the idea of the proof of the estimate from above. It consists of finding an appropriate sequence of
``good'' partitions of $D$. This is done in the following way:
\begin{enumerate}

\item Divide $D$ into a number $m_N^d$ (which is small in comparison with $N$) of equal subregions $D_i^N$. 
On each $D_i^N$, instead of $f$, consider the quadratic part of its Taylor polynomial taken at the center ${\bf{h}}_i$ of
$D_i^N$ (call it $P_2^i(f, {\bf x})$). The error of this intermediate approximation is given in Lemma \ref{L_T_d}.

The choice of $m_N^d$ is governed by the following two reasons: there should be few original subregions in comparison with $N$, but their size should be small enough to provide the small enough error of intermediate approximation of $f$ by  $P_2^i(f, {\bf x})$.

\item We find the parameters of the appropriate partition of $D_i^N$ by minimizing the error of multilinear interpolation of  $P_2^i(f, {\bf x})$ on $D_i^N$, $i=1,\dots, m_N^d$.
We choose the size of $R$ in such a way that the overall number $n_i^N$ of elements of partition used for
$D_i^N$ is such that the sum $\displaystyle \sum_{i=1}^{m_N^d}n_i^N$ is approximately $N$, and the errors of
interpolation on each $D_i^N$ are approximately equal.

\item The final partition of $D$ is the union of partitions of each region $D_i^N$, $i=1,\dots, m_N^d$.

\item We show that the sequence of partitions which is optimal for the intermediate approximant (piecewise quadratic function $P_2(f, {\bf x})$) will
be asymptotically optimal for the original function $f$.

\end{enumerate}

Having constructed a partition for the fixed $N$, we interpolate $f$ at the vertices of this partition. 
This will produce a multilinear spline which will be discontinuous along ``small'' number of edges in the partition. Repeating the construction for every $N$, we shall obtain a sequence of partitions and therefore a sequence of interpolating multilinear splines which will be asymptotically optimal.

If we ``give up'' the interpolation at some points (``small'' amount of them) of the sequence of partitions for the sake of having a sequence of continuous multilinear splines, then we shall refine the obtained on each step partition and then ``glue'' splines on the neighboring elements together. The resulting continuous spline we shall call {\it quasi interpolating} spline (see Section \ref{5.2} for detailed construction).

\section{Auxiliary statements.} \label{aux}

The proofs of the following auxiliary statements are straightforward. Similar statements have been proved, for instance, in ~\cite{PhD, us}.

Set for $f \in L_{\infty}(D)$
, $D \subset \RR^d$, 
\begin{equation} \label{mmodulus}
\omega(f,\delta):=\sup \{|f({\xx})-f({\xx}')|: \;\; |{\xx}-{\xx}'|\leq \delta,\;\; {\bf x}, {\bf x'} \in D  \}, \;\; \delta \geq 0,
\end{equation}
where $|{\bf x}|:=\displaystyle \max_{1\leq i\leq d} |x_i|$ for ${\bf x} \in \RR^d$.
Set for $f \in C^2(D)$
\begin{equation} \label{mm1}
\omega^*(f,\delta):=\max_{1\leq i,j\leq d}\{\omega(f_{x_ix_j}, \delta)\},
\end{equation}
where $f_{x_ix_j}$ denotes mixed derivative of $f$ with respect to variables $x_j$ and $x_i$, $i,j=1,\dots,d$.

\begin{lemma} \label{L_T_d}
Let $f \in C^2(D)$ and $P_2({\bf x})$ denote its quadratic Taylor polynomial at the center ${\bf x}_0$ of a
cube $D_h\subset D$ in $\RR^d$ with side length equal to $h$. Then we have the following estimate:
\begin{equation} \label{eest}
|f({\bf x})-P_2({\bf x})| \leq \frac{d^2}{2} \left(\frac{h}{2}\right)^2 \omega^* \left (f,\frac h2\right ),\;\;\;\;\; {\bf
x} \in D_h,
\end{equation}
where $\omega^*(f,t)$ is defined in (\ref{mm1}).
\end{lemma}

For a fixed ${\bf a} \in \RR^d$ and an arbitrary box $R$ denote
$$
R+{\bf a}=\{{\bf x}+{\bf a}, \;\;\; {\bf x}\in R \}.
$$

\begin{lemma} \label{invar_b}
For the given quadratic function
\begin{equation} \label{f_L}
Q({\bf x})=\displaystyle \sum_{i=1}^d A_i x_i^2,
\end{equation}
any box $R$, and ${\bf a}\in \RR^d $ errors (in any $L_p$ norm) of multilinear interpolation of $Q({\bf x})$ at the vertices of $R$ and $R+{\bf a}$ are equal.
\end{lemma}

\begin{lemma} \label{interp}
The interpolant of the quadratic function (\ref{f_L})
on the $d$-dimensional box $R^d:=\displaystyle \prod_{i=1}^d [-h_i,h_i]$ is a constant function
$$
s(Q,R^d):=\displaystyle \sum_{i=1}^d A_i h_i^2.
$$
\end{lemma}
Proofs of last two statements are simple linear algebra exercises.

\section{Interpolation of quadratic functions with arbitrary signature.} \label{quadr}

Let the quadratic form
\begin{equation} \label{b-q}
Q({\bf x})=\displaystyle \sum_{i=1}^k x_i^2-\displaystyle \sum_{i=k+1}^d x_i^2
\end{equation}
for $0 \leq k \leq d$ be given. 

\begin{lemma} \label{arb_sign}
Let $Q$ be the quadratic function of form (\ref{b-q}),
let $R^d:=\displaystyle \prod_{i=1}^d [-h_i,h_i]$, and let $P$ be the unique  polynomial from ${\cal P}^*_1$ interpolating $Q$ at the vertices of $R^d$. Then
\begin{equation} \label{sum}
\|Q-P\|_{L_{\infty}(R^d)}=\displaystyle \max \left \{ h_1^2+\dots+h_k^2, h_{k+1}^2+\dots+h_d^2 \right \}.
\end{equation}
\end{lemma}

{\bf {Proof}:} We shall proceed by induction on the number of variables. To prove the basis of induction we need to
consider two cases $d=2$, $k=1$ and $d=2$,  $k=2$. 

{\underline {Case 1}.} Let
the quadratic form
$$
Q(x,y)=x^2-y^2
$$
 and an arbitrary rectangle $R^2:=[-h_1,h_1]\times [-h_2,h_2]$ be given.

By Lemma ~\ref{interp}, the bilinear interpolant to the function $Q(x,y)$ on the rectangle $R$, denoted by $T_{Q,R}(x,y)$, is a
constant equal to
$$
T_{Q,R}(x,y)=h_1^2-h_2^2.
$$
Observe that, because of the symmetry, the error in the uniform norm on $R$ is the same as the error on
$[0,h_2]\times[0,h_1]$.

Denote the difference between the function $Q(x,y)$ and the interpolant $T_{Q,R}(x,y)$ by
\begin{equation} \label{delta}
\delta(x,y):=x^2-y^2-h_1^2+h_2^2.
\end{equation}
Clearly, the point $(0,0)$ is the only critical critical point in $[0,h_2]\times[0,h_1]$ of this function. The value of the difference (\ref{delta}) at this
point is
$$
|\delta(0,0)|=|h_1^2-h_2^2|.
$$
In addition, observe that on the boundary of $[0,h_2]\times[0,h_1]$ we have
$$
\delta(x,h_2)=x^2-h_1^2,\;\;\;\hbox{and}\;\;\;\delta(h_1,y)=y^2-h_2^2
$$
and, hence, maximal values of the difference are
$$
|\delta(0,h_2)|=h_1^2 \;\;\;\hbox{and}\;\;\;|\delta(h_1,0)|=h_2^2.
$$
Therefore, the error in the uniform norm can be rewritten as follows:
\begin{equation} \label{min_1}
\displaystyle \displaystyle \max\{|h_1^2-h_2^2|, h_1^2, h_2^2 \}=\displaystyle \max\{h_1^2, h_2^2 \},
\end{equation}
and the statement of the Lemma is proved in the Case 1.

{\underline {Case 2}.} Similarly to the Case 1 we can show that the error of interpolation of
the quadratic form
$$
Q(x,y)=x^2+y^2
$$
on an arbitrary rectangle $R=[-h_1,h_1]\times [-h_2,h_2]$ by multilinear polynomial which in this case is going to have a form
$$
T_{Q,R}(x,y)=h_1^2+h_2^2.
$$
In addition, note that the error in the uniform norm on the rectangle $R$ is the same as the error on $[0,h_2]\times[0,h_1]$.

Denote the difference between function $Q(x,y)$ and interpolant $T_{Q,R}(x,y)$ by
\begin{equation} \label{delta_dif}
\delta(x,y):=x^2+y^2-h_1^2-h_2^2.
\end{equation}
Clearly, the point $(0,0)$ is the only critical point of this function inside $R$. The value of  difference (\ref{delta_dif}) at this point is
$$
\delta(0,0)=h_1^2+h_2^2.
$$
In addition, observe that on the boundary of $[0,h_2]\times[0,h_1]$ we have
$$
\delta(x,h_2)=x^2-h_1^2,\;\;\;\hbox{and}\;\;\;\delta(h_1,y)=y^2-h_2^2
$$
and, hence, the maximal values are
$$
|\delta(0,h_2)|=h_1^2 \;\;\;\hbox{and}\;\;\;|\delta(h_1,0)|=h_2^2.
$$
Therefore,
in  Case 2 the error 
is equal to
\begin{equation} \label{min_1_2}
\displaystyle \max\{h_1^2+h_2^2, h_1^2, h_2^2 \} = h_1^2+h_2^2.
\end{equation}
This completes the proof of the basis of induction.

Next let us consider form (\ref{b-q}) with signature $(k,d-k)$ when $0<k<d$.  As before we can see that the interpolant to (\ref{b-q}) on $R$ is the constant 
$$
T_{Q,R}({\bf x})=\displaystyle \sum_{i=1}^k h_i^2-\displaystyle \sum_{i=k+1}^d h_i^2.
$$
Denote the difference between the function and the interpolant by
$$
\delta({\bf x}):=Q({\bf x})-T_{Q,R}({\bf x}).
$$
Let us investigate critical points of this function. 
For brevity of computation, let us denote 
$$
S_k:=\displaystyle \sum_{i=1}^k h_i^2, \qquad S_k^{(i)}:=\displaystyle \sum_{j=1, j\neq i}^k h_j^2,\qquad
S_d:=\displaystyle \sum_{i=k+1}^d h_i^2, \qquad S_d^{(i)}:= \displaystyle \sum_{j=k+1, j\neq i}^d h_j^2.
$$
As before, the center ${\bf 0}$ is the only critical point of $\delta({\bf x})$ inside $R$.
The value of the difference at the center is
\begin{equation}
\delta({\bf 0})=\left | S_k-S_d\right|.
\end{equation}
Let us consider the error on the boundary.

On the face $x_i=h_i$ in the case when $i\leq k$ form (\ref{b-q}) becomes
\begin{equation}
\displaystyle \sum_{j=1, j\neq i}^k x_j^2-\displaystyle \sum_{j=k+1}^d x_j^2
\end{equation}
and the error by hypothesis of induction is
\begin{equation}
\max \left \{ S_k^{(i)}, S_d \right \}.
\end{equation}
Similarly, on the face $x_i=h_i$ in the case when $i> k$ form (\ref{b-q}) becomes
\begin{equation}
\displaystyle \sum_{j=1}^k x_j^2-\displaystyle \sum_{j=k+1, j\neq i}^d x_j^2,
\end{equation}
and the error by hypothesis of induction is
\begin{equation}
\max \left \{ S_k, S_d^{(i)} \right \}.
\end{equation}
Therefore, the global error is 
\begin{eqnarray}
\Delta &=& \max \left \{\left | S_k-S_d\right|, \displaystyle \max_{1\leq i\leq k}\max \left \{ S_k^{(i)}, S_d \right \},  \displaystyle \max_{k+1\leq i\leq d}\max \left \{ S_k, S_d^{(i)} \right \}   \right \} \nonumber \\
{} &=& \max \left \{\displaystyle \max_{1\leq i\leq k}\max \left \{ S_k^{(i)}, S_d \right \},  \displaystyle \max_{k+1\leq i \leq d}\max \left \{ S_k, S_d^{(i)} \right \}   \right \} \nonumber \\
{} &=& \max \left \{\max \left \{ \displaystyle \max_{1\leq i\leq k}S_k^{(i)}, S_d \right \},  \max \left \{ S_k, \displaystyle \max_{k+1\leq i \leq d} S_d^{(i)} \right \}   \right \} \nonumber \\
{} &=& \max \left \{\displaystyle \max_{1\leq i\leq k}S_k^{(i)},  S_k, S_d,  \displaystyle \max_{k+1\leq i \leq d} S_d^{(i)}   \right \} \nonumber \\
{} &=& \max \left \{   S_k, S_d \right \}. \nonumber
\end{eqnarray}
The lemma is proved. $\square$

Next we shall compute the minimal value of the error $\Delta$ for the quadratic form (\ref{b-q}) with signature $(k,d-k)$. Denote by
\begin{equation} \label{D_min}
\tilde \Delta:=\displaystyle \min_{h_i}  \left \{\displaystyle \sum_{j=1}^k h_j^2, \displaystyle \sum_{j=k+1}^d h_j^2   \right
\},
\end{equation}
where min is taken over all $h_i$ such that the volume
$
2^d \displaystyle \prod_{i=1}^d h_i 
$
is fixed ($=V$).
In the next two lemmas we shall provide the value of $\tilde \Delta$ in the case of the quadratic form (\ref{b-q}) with signature $(k, d-k)$ when $k\neq 0, d$ (Lemma \ref{42})  and in the case of positive definite quadratic form $Q({\bf {x}})=\displaystyle \sum_{i=1}^d x_i^2$ (Lemma
\ref{L4}).

\begin{lemma}\label{42}
The minimal $L_{\infty}$-error of interpolation of quadratic form (\ref{b-q}) by a polynomial from ${\cal P}^*_1$  on all $d$-dimensional boxes of fixed volume $V$ is 
\begin{equation}
\tilde \Delta=\frac 14 k^{\frac kd } (d-k)^{1-\frac kd} V^{\frac 2d}.
\end{equation}
\end{lemma}

{\bf {Proof}:} 
Due to Lemma ~\ref{invar_b} we may consider only boxes $R^d=\displaystyle \prod_{i=1}^d [-h_i,h_i]$ centered at the origin.

To minimize the $l_{\infty}$ -norm of the vector $\left ( \displaystyle \sum_{j=1}^k h_j^2, \displaystyle \sum_{j=k+1}^d h_j^2 \right ) $ given by expression in (\ref{sum}), we minimize the $l_p$-norm of this vector
(with an arbitrary $p$) of (\ref{sum}), i.e. the expression
\begin{equation} \label{fp}
\left (\left(\displaystyle \sum_{j=1}^k h_j^2 \right )^p+\left( \displaystyle \sum_{j=k+1}^d h_j^2\right ) ^p \right)^{\frac 1p},
\end{equation}
under assumption that the volume of the box is fixed ($=V$) and take the value
of the minimum when $p=\infty$.

Indeed, if for some set $M\subset \RR^d$ we denote
$$
\|x^*_p\|_{l_p}=\displaystyle \inf_{x\in M} \|x\|_{l_p}=A_p,
$$
$$
\|x^*\|_{l_{\infty}}=\displaystyle \inf_{x\in M} \|x\|_{l_{\infty}}=A_{\infty},
$$
then it is not difficult to check that
\begin{equation}\label{Ap}
\displaystyle \lim_{p\to \infty} A_p=A_{\infty}.
\end{equation}
Indeed, obviously $A_{\infty}\leq A_p$, for any $p$. On the other hand,
$$
A_p=\|x^*_p\|_{l_p}\leq \|x^*\|_{l_p}\leq \|x^*\|_{l_{\infty}}+\varepsilon_p.
$$
Last two estimates combined imply (\ref{Ap}).

%

Denote the length of the sides of the box by $h_i$, $i=1,\dots, d$. The assumption of volume being fixed is equivalent to
\begin{equation} \label{as_volume}
2^{2d}\displaystyle \prod_{i=1}^d h_i^2=V^2.
\end{equation}

The standard routine of minimizing (~\ref{fp}) leads to the only solution of the minimization problem
$$
h_1^2=\dots=h_k^2=:x,
\qquad h_{k+1}^2=\dots=h_d^2=:y.
$$
Taking this into consideration together with assumption (\ref{as_volume}) which now can be rewritten as
$$
x^k y^{d-k}= V^2 2^{-2d},
$$
we can find $x$ and $y$:
\begin{equation}\label{x}
x=\left ( \frac{d-k}{k}\right)^{(1-\frac kd)(1-\frac 1p)} \frac{V^{\frac 2d}}{4}, \qquad \hbox{and} \qquad y=\left (
\frac{d-k}{k}\right)^{-\frac kd(1-\frac 1p)} \frac{V^{\frac 2d}}{4}.
\end{equation}
In the case $p=\infty$ we have
$$
x=\frac 14 \left ( \frac{d-k}{k}\right)^{(1-\frac kd)} V^{\frac 2d}, \qquad \hbox{and} \qquad y=\frac 14 \left (
\frac{d-k}{k}\right)^{-\frac kd} V^{\frac 2d}.
$$
Therefore,
\begin{equation}
h_i=\frac 12 \left ( \frac{d-k}{k}\right)^{\frac{d-k}{2d}} V^{\frac 1d}, \;\;\; i\leq k,
\end{equation}
\begin{equation}
h_j=\frac 12 \left ( \frac{d-k}{k}\right)^{-\frac {k}{2d}} V^{\frac 1d},\;\;\; j>k.
\end{equation}
Hence, the minimal value $\tilde \Delta$ of the error $\Delta$ is
\begin{equation}
\tilde \Delta=\frac 14 k\left ( \frac{d-k}{k}\right)^{(1-\frac kd)} V^{\frac 2d}=\frac 14 k^{\frac kd } (d-k)^{1-\frac kd}
V^{\frac 2d}. \;\;\;
\end{equation}
$\square$

\begin{lemma} \label{L4}
The minimal $L_{\infty}$-error of interpolation of positive definite quadratic form by a polynomial from ${\cal P}_1$  on all $d$-dimensional boxes of fixed volume $V$ is 
\begin{equation}
\tilde \Delta=\frac d4  V^{\frac 2d}.
\end{equation}
\end{lemma}

{\bf {Proof}:} 
 Clearly, the minimum of the function $\displaystyle \sum_{i=1}^d h_i^2$ with additional assumption
(\ref{as_volume}) is achieved when all $h_i$ are equal, i.e.
$
h_1=h_2=\dots=h_d:=h.
$
From condition (\ref{as_volume}) we also have
$$
h=\frac{V^{\frac 1d}}{2}, \;\; \hbox{and, hence,}\;\; \tilde \Delta=\displaystyle \min_{h_i}\left \{ \displaystyle \sum_{i=1}^d h_i^2 \right \}=d h^2=d \frac{V^{\frac 2d}}{4}. $$
$\square$

Now let the quadratic form
\begin{equation} \label{b-qA}
Q({\bf x})=\displaystyle \sum_{i=1}^d A_i x_i^2
\end{equation}
with $A_i>0$ for all $0\leq i\leq k$ and $A_i<0$ for all $k+1\leq i\leq d$ be given. 

\begin{lemma} \label{arb_sign_A}
The $L_{\infty}$ - error of interpolation of quadratic form (\ref{b-qA}) by polynomials ${\cal P}_1$ on
the $d$-dimensional box $P$ of volume $V(P)$ is
\begin{equation} \label{b_er}
\frac 14 k^{\frac kd } (d-k)^{1-\frac kd} \left (V(P) \sqrt{\displaystyle \prod_{i=1}^d |A_i|}\right)^{\frac 2d}.
\end{equation}
\end{lemma}

{\bf {Proof}:} For the given quadratic form $Q({\bf x})=\displaystyle \sum_{i=1}^d |A_i| x_i^2$ let us consider a linear
transformation $F$ such that
\begin{equation} \label{b_trans}
(Q\circ F)({\bf u})=\displaystyle \sum_{i=1}^d u_i^2.
\end{equation}
In other words,
\begin{equation}\label{b_another}
F({\bf u})=\left (\frac{u_1}{\sqrt{|A_1|}},\dots, \frac{u_d}{\sqrt{|A_d|}}  \right ).
\end{equation}
Observe that the determinant of the inverse of this transformation is
\begin{equation} \label{det}
\det (F^{-1})=\sqrt{\displaystyle \prod_{i=1}^d |A_i|}.
\end{equation}

Let us consider the box $F^{-1}(P)$ which clearly has the volume
\begin{equation} \label{b_volume} V(F^{-1}(P))=V(P) \det
(F^{-1}).
\end{equation}
Combining the result of the previous lemma about the error of interpolation on the box $F^{-1}(P)$ with
(\ref{b_volume}) and (\ref{det}), we obtain (\ref{b_er}). $\square$

Similarly, in the case of positive definite form we obtain the following statement.

\begin{lemma} \label{pos_sign_A}
The error of interpolation of the positive definite quadratic form  by polynomials ${\cal P}_1$ on the $d$-dimensional box
$P$ of volume $V(P)$ is equal to
\begin{equation} \label{b_er1}
\frac d4  \left (V (P) \sqrt{\displaystyle \prod_{i=1}^d A_i}\right)^{\frac 2d}.
\end{equation}
\end{lemma}

\section{Error of interpolation of $C^2$ functions defined on $[0,1]^d$. Estimate from above.} \label{sec_bil}
\subsection{Estimate from above for interpolating splines.} \label{sec_bil_1}
\begin{lemma}
Let $f \in C^2_k(D)$. Then
\begin{equation} \label{ab}
\displaystyle \limsup_{N \to \infty} \frac{N^{\frac 2d}{\mathbf R}_N(f)}{\gamma_{k,d}\left (\displaystyle \int_D |H(f;{\bf{ x}})|^{\frac 12}\Omega({\bf x})^{\frac d2} d{\bf x}\right ) ^{\frac 2d}}\leq 1.
\end{equation}
\end{lemma}

{\bf Proof.}

 Let $f\in C^2(D)$ be given. For a fixed $\varepsilon \in (0,1)$ and for every $N$ we define
\begin{equation} \label{dmm_N_b}
m_N:=\min\left \{ m>0: \;\; \frac{d^2}{2} \left(\frac{1}{2m}\right)^{2} \omega ^* \left (f, \frac {1}{2m} \right ) \leq
\frac{\varepsilon}{N^{\frac 2d}} \right \},
\end{equation}
where $\omega^*(f, \delta)$ is a function defined in (\ref{mm1}). Observe that for $m_N$ defined in such a way it
is true that $m_N \to \infty$ as $N \to \infty$. In addition,
\begin{equation} \label{dsstar}
\frac{N^{\frac 2d}}{m_N^{2}} \to \infty, \;\;\; \hbox{as} \;\;\; N \to \infty,
\end{equation}
i.e. $m_N=o(N^{\frac{1}{{d}}})$ as $N\to \infty$. Indeed, by the definition of $m_N$ for all large enough $N$ we have
$$
\frac{N^{\frac 2d}}{m_N^2}= \frac{(m_N-1)^2}{m_N^2} \frac{N^{\frac 2d}}{(m_N-1)^2}\geq \varepsilon \frac{8}{d^2}\frac{(m_N-1)^2}{m_N^2} \left(\omega^*\left(f, \frac{1}{2(m_N-1)}\right )\right)^{-1} \to \infty, \;\;
\hbox{as}\;\; N\to \infty,
$$
since $\left(\frac{m_N-1}{m_N}\right)^2 \to 1$ and $\omega^*\left(f, \frac{1}{2(m_N-1)}\right) \to 0$ as $ N\to \infty$.
Hence, (\ref{dsstar}) is proved.

 Let us divide the unit cube $D$ into cubes with side length equal to
$\frac{1}{m_N}$ and denote the resulting cubes by $D_l^N$, $l=1,\dots,m_N^d$. Next let us take the center point ${\bf
x}_l^N$ in each cube $D_l^N$  and set
$$
A_{i,j}^{N,l}:=\frac 12 \frac{\partial^2 f}{\partial x_i \partial x_j}({\bf x}_l^N),\;\;\; i,j=1,\dots,d,\;\;\;
l=1,\dots,m_N^d.
$$
In addition, denote by
$$
H({\bf x}_l^N):=\displaystyle \prod_{i=1}^d A_{i,i}^{N,l},\;\;\; \;\;\; l=1,\dots,m_N^d.
$$
Set the number of elements to be used on $D_l^N$ to be
\begin{equation} \label{dnni}
n_l^N:=\left [ \frac{N(1-\varepsilon)|H({\bf x}_l^N)|^{\frac{1}{2}}\Omega({\bf{x}}_l^N)^{\frac d2}} {\displaystyle
\sum_{j=1}^{m_N^d} |H({\bf x}_j^N)|^{\frac{1}{2}}\Omega({\bf{x}}_j^N)^{\frac d2}}\right],\;\;\;l=1,\dots,m_N^d.
\end{equation}

It is essential that $n_l^N \to \infty$ when $N \to \infty$. This follows from the estimate
\begin{equation} \label{ni_est}
n_l^N \geq \left [ \frac{N(1-\varepsilon) {\displaystyle \min_{{\bf x} \in D} \{ |H(f; {\bf x}})|\}^{\frac{1}{2}} \displaystyle \min_{{\bf x} \in D} \{ \Omega({\bf x})
\}^{\frac{d}{2}}} {m_N^d\|H\|_{\infty}^{\frac{1}{2}}\|\Omega\|_{\infty}^{\frac{d}{2}}}\right ],
\end{equation}
together with (\ref{dsstar}), and the fact that $\displaystyle \min_{{\bf x} \in D} \{ \Omega({\bf x})\}>0$ and $\displaystyle \min_{{\bf x} \in D} \{ |H(f; {\bf x})|\}>0$.

Now for $k\neq 0, d$ let us set 
\begin{equation} \label{hi1}
h_{l}^N:=\frac 12 \left ( \frac{d-k}{k}\right)^{\frac{d-k}{2d}}\left (\frac{1}{m_N^d n_l^N}\right)^{\frac 1d}, \;\;\;
1\leq i\leq k,\;\; l=1,\dots, m_N^d,
\end{equation}
\begin{equation}\label{hi2}
\tilde h_{l}^N:=\frac 12 \left ( \frac{d-k}{k}\right)^{-\frac {k}{2d}} \left (\frac{1}{m_N^d n_l^N}\right)^{\frac 1d},\;\;\;
k+1\leq i \leq d, \;\; l=1,\dots, m_N^d.
\end{equation}
For the positive definite form, i.e. $k=0$ or $k=d$ we set
\begin{equation}\label{hi3}
h_{l}^N=\tilde h_{l}^N:=\frac 12 \left (\frac{1}{m_N^d n_l^N}\right)^{\frac 1d}, \;\;\;  l=1,\dots, m_N^d
\end{equation}
The intersection of the lattice
\begin{equation}
[Lh_{l}^N,(L+1)h_{l}^N]^k\times [L\tilde h_{l}^N,(L+1)\tilde h_{l}^N]^{d-k}, \;\;\; L \in \ZZ,
\;\;l=1,\dots, m_N^d,
\end{equation}
with $D_l^N$ provides the partition of $D_l^N$, $l=1,...,m_N^d$.

The union of partitions of each $D_l^N$, $l=1,...,m_N^d$, provide the partition of $D$.
Denote it by $\square^*_N(D)$. 
It can be easily seen that due to definitions (\ref{hi1}) and (\ref{hi2}), combined with the estimate (\ref{ni_est}), the constructed partition satisfies (\ref{asss}) and, hence, is admissible.

Let us show that the sequence of obtained in such a way
partitions $\{\square^*_N(D)\}_{N=1}^{\infty}$ will be asymptotically optimal.

By $f_N$ denote the piecewise quadratic function constructed in the following way. On $D_1^N$ we set $f_N$ to be
$\displaystyle \sum _{i=1}^{d}A_{i,i}^{N,1}x_i^2$. Then for $l
>1$ on $D_l^N\setminus \displaystyle \cup_{j=1}^{l-1}D_j^N$ we set
$$
f_N({\bf x}):=\displaystyle \sum _{i=1}^{d}A_{i,i}^{N,l}x_i^2.
$$

To estimate the error ${\mathbf R}_N(f)$ we observe that
$$
{\mathbf R}_N(f)\leq \|f - s(f, \square^*_N)\|_{\infty,\Omega} \leq \|f-f_N\|_{\infty,\Omega}+\|f_N-s(f_N, \square^*_N)\|_{\infty,\Omega}
 $$
 $$
 +\|s(f_N,\square^*_N)-s(f, \square^*_N)\|_{\infty,\Omega}\leq 2 \|f-f_N\|_{\infty,\Omega}+\|f_N-s(f_N, \square^*_N)\|_{\infty,\Omega}.
$$
Let us estimate each term. First of all, by Lemma \ref{L_T_d} and the definition of $m_N$ we have
$$
\|f-f_N\|_{\infty,\Omega} \leq \frac{d^2}{2} \left(\frac{1}{2m_N}\right)^{2} \omega^*\left(f, \frac{1}{2m_N}\right)
\|\Omega\|_{\infty}\le \frac{\varepsilon}{N^{\frac 2d}}\|\Omega\|_{\infty}.
$$
Let us estimate the second term now.
Let  $R_l^N \in \square^*_N(D_l^N)$  be an arbitrary element.
By Lemma \ref{arb_sign_A}, for every ${\bf{x}} \in R_l^N$ we have
$$
|f_N({\bf{x}})-s(f_N, \square^*_N; {\bf{x}})|\Omega({\bf x}_l^N)\leq \|f_N-s(f_N, \square^*_N; \cdot)\|_{L_{\infty}(R_l^N)} \Omega({\bf x}_l^N)
$$
$$
=\gamma_{k,d} \left (\frac{1}{m_N^dn_l^N}\sqrt{|H({\bf{x}}_l^N)|}\right)^{\frac
2d}\Omega({\bf x}_l^N).
$$
By the definition of $n_l^N$, for all large enough $N$, for all $l$, and for all ${\bf x} \in R^N_l$, we have
$$
|f_N({\bf x})-s(f_N, \square^*_N; {\bf x})|\Omega({\bf x}_l^N) \leq
$$
$$
\leq \gamma_{k,d} \left (\frac{\displaystyle \sum_{j=1}^{m_N^d} {|H({\bf
x}_j^N)|}^{\frac{1}{2}}\Omega({\bf x}_j^N)^{\frac d2}} {m_N^d N(1-\varepsilon){|H({\bf
x}_l^N)|}^{\frac{1}{2}}\Omega({\bf x}_l^N)^{\frac d2}}\sqrt{|H({\bf{x}}_l^N)|} \right)^{\frac 2d}\Omega({\bf x}_l^N).
$$
Since this estimate does not depend on ${\bf x}$, we obtain
$$
\|f_N-s(f_N, \square^*_N)\|_{\infty, \Omega} \leq \frac {\gamma_{k,d} }{( N(1-\varepsilon))^{2/d}}
\left (\frac{1}{m_N^2}{\displaystyle \sum_{j=1}^{m_N^d} {|H({\bf x}_j^N)|}^{\frac{1}{2}}\Omega({\bf x}_j^N)^{\frac d2}}
\right )^{\frac 2d}.
$$
Note that
$$
\frac{1}{m_N^2}{\displaystyle \sum_{j=1}^{m_N^d} {|H({\bf x}_j^N)|}^{\frac{1}{2}}\Omega({\bf x}_j^N)^{\frac d2}} \to
\displaystyle \int _D {|H(f;{\bf{x}})|^{\frac 12}\Omega({\bf x})^{\frac d2}}d{\bf x},\;\;\hbox{as}\;\;N\to\infty.
$$
Hence, for all $N$ large enough we have
$$
\|f_N-s(f_N, \square^*_N)\|_{\infty, \Omega} \leq \frac {\gamma_{k,d} }{( N(1-\varepsilon))^{2/d}}
\left ( \displaystyle \int _D {|H(f;{\bf{x}})|^{\frac 12}\Omega({\bf x})^{\frac d2}}d{\bf x} \right)^{\frac 2d}.
$$
Therefore,
$$
\|f-s(f, \square^*_N)\|_{\infty, \Omega} \leq \frac {2\varepsilon}{N^{\frac 2d}}\|\Omega\|_{\infty}+ \frac {\gamma_{k,d} }{( N(1-\varepsilon))^{{\frac 2d}}} \left ( \displaystyle \int _D {|H(f;{\bf{x}})|^{\frac
12}\Omega({\bf x})^{\frac d2}}d{\bf x} \right)^{\frac 2d}.
$$
Because $\varepsilon >0$ is arbitrary,  we obtain the desired estimate from above (\ref{ab}) for ${\mathbf R}_N(f)$.

\subsection{Construction of asymptotically optimal sequence of continuous quasi interpolating splines. }\label{5.2}

In this section, we shall refine the sequence of partitions $\{\square^*_N\}_{N=1}^{\infty}$ constructed following the algorithm in the previous section, to obtain a sequence of partitions and a sequence of continuous splines $\{\tilde s(f, \square_N^*)\}_{N=1}^{\infty}$ on the partitions which will interpolate $f$ at all but $o(N)$ as $N \to \infty$ points.

Let parameters of the grid $h_{l}^N$ and $\tilde h_{l}^N$ be as defined in (\ref{hi1}) and (\ref{hi2}) for $k\neq 0,d$, and in (\ref{hi3}) for $k=0,d$, respectively. Recall that the intersection of the lattice
\begin{equation}
[Lh_{l}^N,(L+1)h_{l}^N]^k\times [L\tilde h_{l}^N,(L+1)\tilde h_{l}^N]^{d-k}, \;\;\; L \in \ZZ,
\;\;l=1,\dots, m_N^d,
\end{equation}
with $D_l^N$ provides the partition of $D_l^N$, $l=1,...,m_N^d$.

\begin{figure}
	\begin{center}
		\includegraphics[scale=0.5]{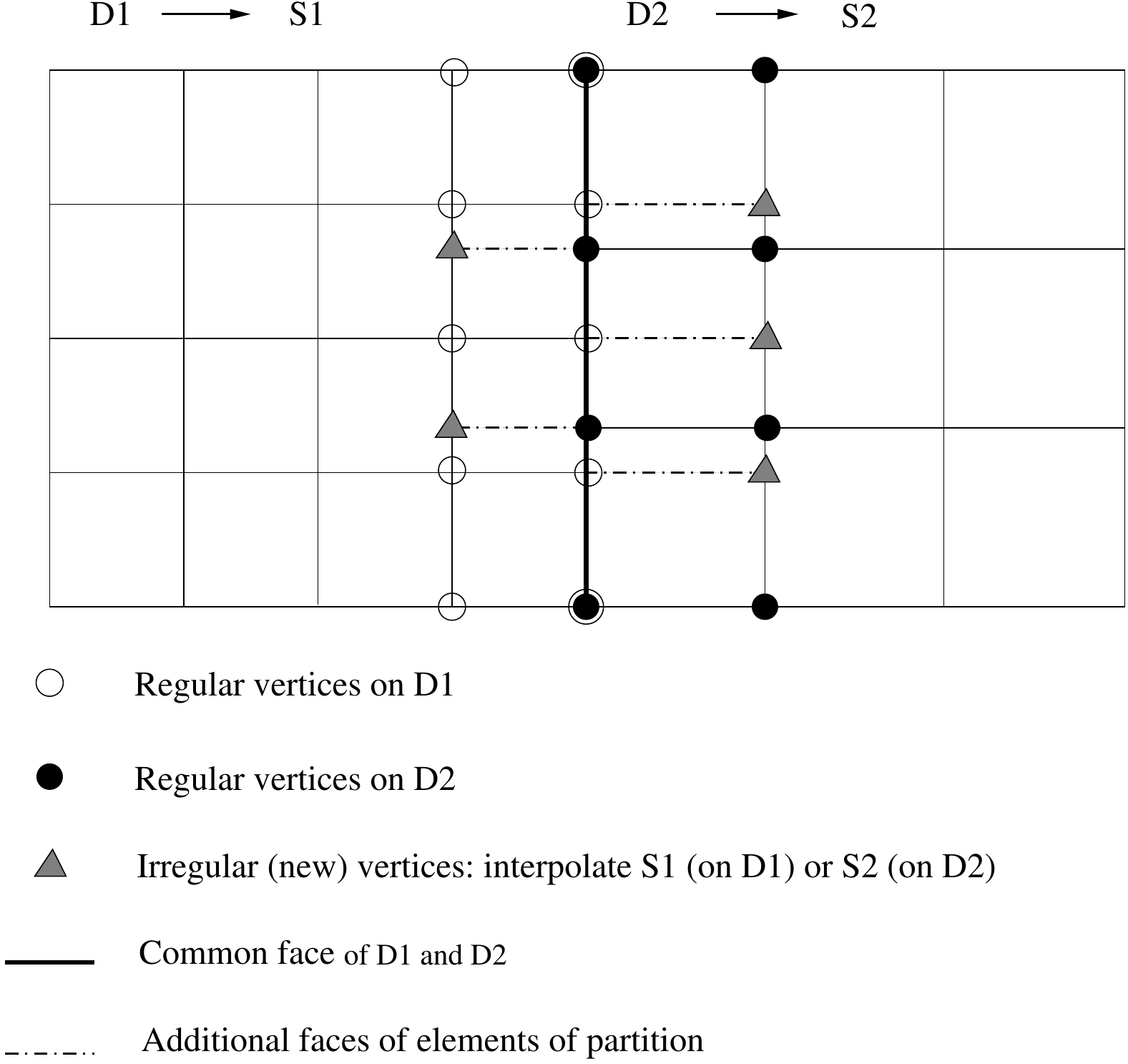}
	\caption{Stitching of partitions on neighboring regions in the case $d=2$.}
	\label{fig:2part}
	\end{center}
	\end{figure}

Let us consider two neighboring $d$-dimensional regions $D^1$ and $D^2$ with corresponding interpolants $S^1$
and $S^2$. The partition of $D^1$ consists of parallelepipeds $\{D_i^1\}$ and the partition of $D^2$ consists of parallelepipeds $\{D_j^2\}$.
 If the parameters of the grid on $D^1$ and $D^2$ are different, we have to subdivide parallelepipeds that have nonempty intersection
  with the common face to ensure the global continuity of the approximant. Let us assume that the common face lies in the $(d-1)$-dimensional coordinate plane.
  
  By $\{D^1_{i,1}\}$ and $\{D^2_{j,1}\}$ let us denote parallelepipeds which have nonempty intersection with the common face of $D^1$ and $D^2$. In addition, by $\{\overline {D^1_{i,1}}\}$ and $\{\overline {D^2_{j,1}}\}$ we shall denote their $(d-1)$-dimensional faces contained in $D^1\cap D^2$. The set of all possible intersections $\overline {D^1_{i,1}}\cap \overline {D^2_{j,1}}$ constitute the partition of $D^1\cap D^2$, which is a refinement of each of partitions $\{\overline {D^1_{i,1}}\}$ and $\{\overline {D^2_{j,1}}\}$.
  
  For each $\overline {D^1_{i,1}}\cap \overline {D^2_{j,1}}$ we consider $D^1_{i,1}\cap \{(\overline {D^1_{i,1}}\cap \overline {D^2_{j,1}})\times \RR\}$ and $D^2_{j,1}\cap \{(\overline {D^1_{i,1}}\cap \overline {D^2_{j,1}})\times \RR\}$, the set of which constitute a refinement of $\{D^1_{i,1}\}$ and $\{D^2_{j,1}\}$, respectively. The new vertices of the refined partition we shall call ``irregular'' to the contrast with vertices of the non refined partitions which we shall refer to as ``regular''.
  
  The continuous spline $S$ on two neighboring elements $D^1$ and $D^2$ is constructed now as follows: $S$ interpolates $S^1$ at all (regular and irregular) vertices of the partition of $D^1$, and  $S$ interpolates $S^2$ at all (regular and irregular) vertices of the partition of $D^2$. Note that automatically the multilinear spline $S$ interpolates $f$, in particular, at all the vertices of the refined partition $\{\overline {D^1_{i,1}}\cap \overline {D^2_{j,1}}\}$ of the common face $D^1\cap D^2$.



The illustration of ``stitching'' of partitions on neighboring regions  in 2-dimensional case is shown in Figure \ref{fig:2part}.



%


Denote the final partition of $D$ by $\tilde \square^*_N(D)$. 
It can be easily seen that due to (\ref{hi1}) and (\ref{hi2}), or (\ref{hi3}), the constructed partition satisfies (\ref{asss}).

Let us compute the number of irregular boxes (denote it by $K_l^N$).  Recalling (\ref{hi1}) and (\ref{hi2}) together
with the fact that the volume of a box from $\square^*_N(D_l^N)$ is $\frac{1}{m_N^dn_l^N}$ we have that
$$
K_l^N \leq c_1 (n_l^N)^{\frac 1d}.
$$
 Hence, the total number of irregular boxes is not
greater than
$$
c_2 \displaystyle \sum _{l=1}^{m_N^d}(n_l^N)^{\frac 1d} \leq c_2 \displaystyle \sum _{l=1}^{m_N^d} \left (
\frac{N(1-\varepsilon)|H({\bf x}_l^N)|^{\frac{1}{2}}\Omega({\bf{x}}_l^N)^{\frac d2}} {\displaystyle \sum_{j=1}^{m_N^d}
|H({\bf x}_j^N)|^{\frac{1}{2}}\Omega({\bf{x}}_j^N)^{\frac d2}} \right )^{\frac 1d}
 \leq c_3
 N^{\frac 1d} {m_N^{d-1}}=o\left( N\right),
$$
as $N \to \infty$ (since $m_N=o\left( N^{\frac 1d}\right )$ as $N\to \infty$ because of (\ref{dsstar})). In particular,
this implies that the number of boxes in the constructed partition will not exceed $N$ for all $N$ large enough.

This verifies the Remark after Theorem 1, because the interpolating spline $s(f_N, \square^*_N)$ constructed in the previous section can be discontinuous only along the faces of the irregular boxes.


It is clear that for two embedded elements of partition the error of multilinear interpolation of a quadratic function is greater on the larger element. Therefore, we shall estimate this error on elements that do not have intersection with the boundary. Hence, the estimate from above for the error of quasi interpolation by splines $\tilde s(f,\tilde \square_N^*)$
can be obtained as in Section \ref{sec_bil_1}.



\section{Error of interpolation of $C^2$ functions defined on $[0,1]^d$. Estimate from below.} \label{6.2}

Let quantities $m_N, n_i^N, D_i^N$ etc. be as defined in the previous section.

\begin{lemma}
Let $f \in C^2_k(D)$. Then
\begin{equation} \label{bel}
\displaystyle \liminf_{N \to \infty} \frac{N^{\frac 2d}{\mathbf R}_N(f)}{\gamma_{k,d}\left (\displaystyle \int_D |H(f;{\bf x})|^{\frac 12}\Omega(\bf x)^{\frac d2} d{\bf x}\right ) ^{\frac 2d}}\geq 1.
\end{equation}
\end{lemma}

{\bf {Proof.}} To obtain the estimate from below we shall consider an arbitrary sequence of admissible partitions, i.e. box
partitions $\{\square_N\}_{N=1}^{\infty}$ which satisfies (\ref{asss}). 

Note that  from (\ref{asss}) it follows that for an arbitrary element $R\in \square_N$, ${\rm{diam}}(R)<\frac{C}{N^{1/d}}$ with some constant $C>0$. Let us consider the $\frac{C}{N^{1/d}}$-neighborhood of the boundary of an arbitrary box $D^N_l$, $l=1,...,m_N^d$. For an arbitrary $\epsilon>0$, the volume of the complement of the $\frac{C}{N^{1/d}}$-neighborhood of the boundary of $D^N_l$, i.e. the ``interior'' of $D^N_l$, is 
\begin{equation}\label{niz}
\left(\frac{1}{m_N}-\frac{2C}{N^{1/d}}\right)^d=\frac{1}{m_N^d}\left(1-\frac{2Cm_N}{N^{1/d}}\right)^d\asymp \frac{1}{m_N^d}\left(1-d\frac{2Cm_N}{N^{1/d}}\right)>\frac{1-\epsilon}{m_N^d}
\end{equation}
for $N$ large enough since $m_N=o(N^{1/d})$ as $N\to \infty$. Therefore, for any $l=1,...,m_N^d$, the sum of volumes of the boxes which have nonempty intersection with the ``interior'' of $D^N_l$ (and, due to (\ref{asss}), lie completely inside of $D^N_l$) is greater than $\frac{1-\epsilon}{m_N^d}$.

Let us show that for any $\epsilon>0$ and for any $N$ large
enough there exists index $l_N$ such that the corresponding  $D_{l_N}^N$ completely contains a ``large enough'' element $R_{l_N}^N \in \square_N$, i.e. element with volume greater than $\frac{1-\epsilon}{m_N^dn_{l_N}^N}$.
Assume to the contrary that there exists $\epsilon_0$ such that for an arbitrary $N_0$ there exists  $N>N_0$ such that for all $i=1,...,m_N^d$ the volume of each box $R\in \square_N$ having nonempty intersection with the ``interior'' of $D^N_i$ (and, therefore, is completely inside of $D^N_i$) is less than or equal to $\frac{1-\epsilon_0}{m_N^dn_i^N}$. 

For each $i=1,...,m_N^d$, by $\nu^N_i$ denote the number of boxes from $\square_N$ that are completely inside of $D^N_i$. Note that
$$
\displaystyle \sum _{i=1}^{m_N^d} \nu^N_i \leq \displaystyle \sum _{i=1}^{m_N^d} n^N_i=N.
$$
This implies that there exists $i^*$ such that $\nu^N_{i^*}\leq n^N_{i^*}$. Hence, taking into consideration the assumption on the volume of each box  that is completely inside of $D^N_{i^*}$ (or have nonempty intersection with the ``interior'' of $D^N_{i^*}$) to be less than or equal to $\frac{1-\epsilon_0}{m_N^dn_i^N}$, the total volume of all such boxes from $\square_N$  is less than or equal to
$$
\frac{1-\epsilon_0}{m_N^dn_i^N}\nu^N_{i^*}\leq \frac{1-\epsilon_0}{m_N^d}
$$
which contradicts to (\ref{niz}).
%


For each $N$ and corresponding $l_N$, set
$$
f_{N, l_N}({\bf x}):=\displaystyle \sum_{i=1}^{d}A_{i,i}^{N,l}x_i^2.
$$
Observe that
$$
\|f-s(f,\square_N)\|_{L_{\infty,\Omega}(R^N_{l_N})}\geq \|f_{N, l_N}-s(f_{N,
l_N},\square_N)\|_{L_{\infty,\Omega}(R^N_{l_N})}-2\|f-f_{N, l_N}\|_{L_{ \infty,\Omega}(R^N_{l_N})}.
$$
By Lemma \ref{arb_sign_A} we have for some $\varepsilon >0$
$$
\|f_{N, l_N}-s(f_{N, i_N},\square_N)\|_{L_{\infty,\Omega}(R^N_{l_N})}\geq (1-\varepsilon)\gamma_{k,d}\left (\frac{1}{m_N^dn_l^N}\sqrt{|H({\bf{x}}_l^N)|}\right)^{\frac 2d}\Omega({\bf x}_l^N).
$$
By the definition of $n_{l_N}^N$  we have that for all $N$ large enough
$$
\begin{array}{lll}
\gamma_{k,d} \left (\frac{1}{m_N^dn_l^N}\sqrt{|H({\bf{x}}_l^N)|}\right)^{\frac
2d}\Omega({\bf x}_l^N)
&=&\gamma_{k,d} \left (\frac{\displaystyle \sum_{j=1}^{m_N^d} {|H({\bf
x}_j^N)|}^{\frac{1}{2}}\Omega({\bf x}_j^N)^{\frac d2}} {m_N^d N(1-\varepsilon){|H({\bf
x}_l^N)|}^{\frac{1}{2}}\Omega({\bf x}_l^N)^{\frac d2}}\sqrt{|H({\bf{x}}_l^N)|} \right)^{\frac 2d}\Omega(x_l^N)\cr
&>&\frac  {\gamma_{k,d}}{ N^{\frac 2d}} \left ( \int_{D}|H(f;{\bf x})|^{\frac12}\Omega({\bf x})^{\frac d2}d{\bf x}\right )^{\frac 2d}.\cr
\end{array}
$$
 Hence, for all $N$ large enough we obtain
$$
\|f_N-s(f_N,\square_N)\|_{\infty,\Omega} > (1-\varepsilon)\frac  {\gamma_{k,d}}{ N^{\frac 2d}}
\left ( \int_{D}|H(f;{\bf x})|^{\frac12}\Omega({\bf x})^{\frac d2}d{\bf x}\right )^{\frac 2d}.
$$
On the other hand
$$
\|f-f_{N, i_N}\|_{L_{\infty,\Omega}(R^N_{i_N})}\leq  \|\Omega\|_{\infty}\frac{\varepsilon}{N^{\frac 2d}}$$ due to the
choice of $m_N$. Hence, we obtain that for all large enough $N$
$$
\|f-s(f,\square_N)\|_{\infty, \Omega} \geq (1-c_4 \varepsilon)\frac  {\gamma_{k,d}}{N^{\frac 2d}}
\left ( \int_{D}|H(f;{\bf x})|^{\frac12}\Omega({\bf x})^{\frac d2}d{\bf x}\right )^{\frac 2d}
$$
with some positive constant $c_4$. Therefore,
$$
\displaystyle   \liminf_{N \to \infty} \|f-s(f,\square_N)\|_{\infty,\Omega}\geq {\frac  {\gamma_{k,d}}{ N^{\frac 2d}} \left ( \int_{D}|H(f;{\bf x})|^{\frac12}\Omega({\bf x})^{\frac d2}d{\bf x}\right )^{\frac 2d}}.\;\; 
$$
$\square$

The estimate from below for quasi interpolating splines can be obtained analogously.





Department of Mathematics and Statistics\\
Sam Houston State University\\
Box 2206\\
Huntsville, TX 77340-2206\\
Phone: 936.294.4884\\
Fax: 936.294.1882\\
Email: babenko@shsu.edu\\


\begin{thebibliography}{99}


\bibitem{Bab}
{V.F. Babenko}, {Interpolation of continuous functions by piecewise linear ones}, Math. Notes, 24, no.1,
(1978) 43--53.


\bibitem{us}
{ V. Babenko, Yu. Babenko, A. Ligun, A. Shumeiko,} {On asymptotical behavior of the optimal linear spline
interpolation error of $C^2$ functions}, East J. Approx., V. 12, N. 1 (2006), 71--101.

\bibitem{PhD}
{Yu. Babenko}, { On the asymptotic behavior of the optimal error of spline interpolation of multivariate functions}, PhD thesis, 2006. 


\bibitem{DeL}
{A. Below, J. De Loera, J. Richter-Gebert},  { The complexity of finding small triangulations of convex
3-polytopes.} SODA 2000 special issue. J. Algorithms 50 (2004), no. 2, 134--167.

\bibitem{Bern}
{M. Bern, D. Eppstein,} { Mesh generation and optimal triangulation}, manuscript.



\bibitem{kodla2}
{M. Bertram, J. Barnes, B. Hamann, K. Joy, H. Pottmann, D. Wushour}, { Piecewise optimal triangulation for
the approximation of scattered data in the plane}, Comput. Aided Geom. Design 17, no. 8, (2000) 767--787.

\bibitem{Bor}
{K. B$\rm{\ddot{o}}$r$\rm{\ddot{o}}$czky}, {Approximation of general smooth convex bodies},  Adv. in Math.,
{153} (2000) 325--341.



\bibitem{KL}
{K. B$\rm{\ddot{o}}$r$\rm{\ddot{o}}$czky, M. Ludwig,}  { Approximation of Convex Bodies and a Momentum Lemma
for Power Diagrams}, Monatshefte f$\it{\ddot{u}}$r Mathematik, V. 127, N. 2, (1999) 101--110.

\bibitem{Chen3}
{L. Chen,} { New analysis of the sphere covering problems and optimal polytope approximation of convex bodies}, J. of Approx. Theory 133 (2005), 134--145.


\bibitem{Chen2}
{L. Chen}, { Optimal Delaunay triangulations}, J. Comp. Math., Vol. 22, No. 2, 2004, 299-308.

\bibitem{Chen1}
{L. Chen, P. Sun, J. Xu} { Optimal anisotropic meshes for minimizing interpolation errors in $L_p$-norm},
Math. Comp.,



\bibitem{Daz1}
{ E. F. D'Azevedo} { Are bilinear quadrilaterals better than linear triangles?} SIAM J. Sci. Comput. 22
(2000), no. 1, 198--217.

\bibitem{Daz2}
{ E. F. D'Azevedo}, { Optimal triangular mesh generation by coordinate transformation.} SIAM J. Sci. Statist.
Comput. 12 (1991), no. 4, 755--786.

\bibitem{Daz3}
{ E. F. D'Azevedo, R. B. Simpson}, {On optimal interpolation triangle incidences}. SIAM J. Sci. Statist.
Comput. 10 (1989), no. 6, 1063--1075.

\bibitem{Delaunay}
{B. N. Delone, S.S. Ryshkov}, {Extremal problems of the theory of positive quadratic forms}. (Russian)
Collection of articles dedicated to Academician Ivan Matveevich Vinogradov on his eightieth birthday, I. Trudy Mat.
Inst. Steklov. 112 (1971), 203--223, 387.


\bibitem{DDI}
{L. Demaret, N. Dyn, A. Iske,} { Image compression by linear splines over adaptive triangulations},  IEEE
Transactions on Image Processing.



\bibitem{Toth}
{L. Fejes Toth,} Lagerungen in der Ebene, auf der Kugel und im Raum, 2nd edn. Berlin: Springer, 1972.



\bibitem{Gr2}
{P. Gruber,} { Volume approximation of convex bodies by inscribed polytopes}, Math. Ann., V. 281, 1988,
pp.229--245.


\bibitem{Gr}
{P. Gruber,}  { Aspects of approximation of convex bodies.} { In: GRUBER P.M., WILLS J (eds) Handbook of
Convex Geometry A, Amsterdam: North-Holland,} (1999) 319--345.

\bibitem {Gr3}
{P. Gruber}, { Error of asymptotic formulae for volume approximation of convex bodies in $E^d$,} Monatsh.
Math. 135 (2002) 279-304.


\bibitem{Huang}
{W. Huang},{ Variational mesh adaptation: isotropy and equidistribution.} J. Comput. Phys. 174 (2001), no. 2,
903--924.

\bibitem{HSun}
{W. Huang, W.  Sun}, { Variational mesh adaptation. II. Error estimates and monitor functions.} J. Comput.
Phys. 184 (2003), no. 2, 619--648



\bibitem{LSh}
{Ligun A.A., Shumeiko A.A.,} Asymptotic methods of curve recovery, {Kiev. Inst. of Math. NAS of Ukraine},
1997. (in Russian)


\bibitem{Nadler}
{E. Nadler}, { Piecewise linear best $L_2$ approximation on triangles,}  in: Chui, C.K., Schumaker, L.L. and
Ward, J.D. (Eds.), Approximation Theory V, Academic Press, (1986) 499--502.


\bibitem{kodla1}
{H. Pottmann, R. Krasauskas, B. Hamann, K. Joy, W. Seibold}, { On piecewise linear approximation of quadratic
functions}, J. Geom. Graph. 4, no. 1, (2000) 31--53.


\bibitem{Rippa}
{S. Rippa}, { Long and thin triangles can be good for linear interpolation}, SIAM J. Num. An., Vol. 29, No.
1 (1992), 257--270.



\end{thebibliography}
\end{document}